\documentclass{article}
 
\usepackage{amsmath,amssymb,amsfonts}

\usepackage{pstcol}

\def\lpro{\slash}
\def\rpro{\backslash}

\newcommand{\lgraft}[2]
{\setlength{\unitlength}{1pt}
\begin{picture}(15,17)
	\put(11,0){#1}
	\put(2,15){#2}
	\put(6,8){$\backslash$}
	\end{picture}}
\newcommand{\rgraft}[2]
{\setlength{\unitlength}{1pt}
\begin{picture}(15,17)
	\put(2,0){#1}
	\put(11,15){#2}
	\put(6,8){$\slash$}
	\end{picture}}

\def\C{\mathbb{C}}
 
\def\A{{\cal{A}}}

\def\H{{\cal{H}}}
\def\Hp{{\cal{H}}^{\gamma}}
\def\He{{\cal{H}}^{e}}
\def\Ha{{\cal{H}}^{\alpha}}
\def\Hat{\widetilde{{\cal{H}}^{\alpha}}}
\def\Hq{{\cal{H}}^{\mathrm{qed}}}

\def\D{\Delta}
\def\Dl{\Delta^{\ltimes}}
\def\DPp{\Delta^p_{\gamma}}
\def\DPe{\Delta^p_e}
\def\Da{\Delta^{\alpha}}
\def\Dat{\widetilde{\Delta^{\alpha}}}
\def\Dp{\Delta^{\gamma}}
\def\De{\Delta^{e}}
\def\Dq{\Delta^{\mathrm{qed}}}

\def\SPe{S^p_e}
\def\Sp{S^{\gamma}}

\def\dd{\delta}
\def\ddl{\delta^{\ltimes}}
\def\dds{\delta^{\sigma}}
\def\ddt{\widetilde{\delta}}
\def\ddp{\delta^{\gamma}}
\def\dde{\delta^e}

\def\Up{U^{\gamma}}
\def\Ue{U^{e}}
\def\Rp{R^{\gamma}}
\def\Re{R^{e}}
\def\Cp{C^{\gamma}}
\def\Ce{C^{e}}

\def\ldot{\cdot_{\ltimes}}
\def\sdot{\cdot_{\sigma}}

\def\<{\langle}
\def\>{\rangle}

\def\Id{\mathrm{Id}}
\def\Hom{\mathrm{Hom}}

\newtheorem{ela}{}[section]
\newenvironment{ale}{\begin{ela} \em}{\end{ela}}
\newenvironment{proof}{\noindent{\em Proof.\/}}{\hfill$\square$\par\vspace{.2cm}}
\newenvironment{proof of}[1]{\noindent{\em Proof of (\ref{#1}).\/}}
	{\hfill$\Box$\par\vspace{.2cm}}
\def\@begintheorem#1#2{\trivlist
    \item[\hskip \labelsep{\bfseries #2\ #1}]%
    \ifnum\@defflag=0\itshape\else\upshape\fi\gdef\@defflag{0}}
\def\@opargbegintheorem#1#2#3{\trivlist
    \item[\hskip \labelsep{\bfseries #2\ #1\ (#3)}]\itshape}
\def\newdef#1#2{\newtheorem{#1}{{#2}\gdef\@defflag{1}}}
\def\@defflag{0}
\newcommand{\formula}[1]%
{\par\vspace{.3cm}\noindent$\ds{#1}$\par\vspace{.3cm}\noindent\ignorespaces}
\makeatletter
\def\@eqnnum{\hb@xt@.01\p@{}%
                      \rlap{\normalfont\normalcolor
                        \hskip -\displaywidth(\theela)}}
\def\eqnarray{%
   \stepcounter{ela}%
   \def\@currentlabel{\p@equation\theela}%
   \global\@eqnswtrue
   \m@th
   \global\@eqcnt\z@
   \tabskip\@centering
  \let\\\@eqncr
   $$\everycr{}\halign to\displaywidth\bgroup
       \hskip\@centering$\displaystyle\tabskip\z@skip{##}$\@eqnsel
      &\global\@eqcnt\@ne\hskip \tw@\arraycolsep \hfil${##}$\hfil
      &\global\@eqcnt\tw@ \hskip \tw@\arraycolsep
         $\displaystyle{##}$\hfil\tabskip\@centering
      &\global\@eqcnt\thr@@ \hb@xt@\z@\bgroup\hss##\egroup
         \tabskip\z@skip
      \cr
}

\def\endeqnarray{%
      \@@eqncr
      \egroup
      \global\advance\c@ela\m@ne
   $$\@ignoretrue
}
\def\@@eqncr{\let\reserved@a\relax
    \ifcase\@eqcnt \def\reserved@a{& & &}\or \def\reserved@a{& &}%
     \or \def\reserved@a{&}\else
       \let\reserved@a\@empty
       \@latex@error{Too many columns in eqnarray environment}\@ehc\fi
     \reserved@a \if@eqnsw\@eqnnum\stepcounter{ela}\fi
     \global\@eqnswtrue\global\@eqcnt\z@\cr}

\makeatother
\renewcommand{\|}{
\setlength{\unitlength}{3pt}
\psset{unit=3pt}
\psset{runit=2pt}
\psset{linewidth=0.2}
\begin{pspicture}(0,.5)(2.5,4)
\psline(1,0)(1,3)
\end{pspicture}}

\newcommand{\Y}{
\setlength{\unitlength}{3pt}
\psset{unit=3pt}
\psset{runit=2pt}
\psset{linewidth=0.2}
\begin{pspicture}(-1,.5)(3,4)
\psline(1,0)(1,2)
\psline(1,2)(0,3)
\psline(1,2)(2,3)
\end{pspicture}}

\newcommand\deuxun{
\setlength{\unitlength}{3pt}
\psset{unit=3pt}
\psset{runit=2pt}
\psset{linewidth=0.2}
\begin{pspicture}(0,.5)(5,5)
\psline(3,0)(3,2)
\psline(3,2)(1,4)
\psline(3,2)(4,3)
\psline(2,3)(3,4)
\end{pspicture}}

\newcommand\deuxdeux{
\setlength{\unitlength}{3pt}
\psset{unit=3pt}
\psset{runit=2pt}
\psset{linewidth=0.2}
\begin{pspicture}(1,.5)(6,5)
\psline(3,0)(3,2)
\psline(3,2)(5,4)
\psline(3,2)(2,3)
\psline(4,3)(3,4)
\end{pspicture}}

\newcommand\troisun{
\setlength{\unitlength}{3pt}
\psset{unit=3pt}
\psset{runit=2pt}
\psset{linewidth=0.2}
\begin{pspicture}(-1,.5)(5,6)
\psline(3,0)(3,2)
\psline(3,2)(0,5)
\psline(3,2)(4,3)
\psline(2,3)(3,4)
\psline(1,4)(2,5)
\end{pspicture}}

\newcommand\troisdeux{
\setlength{\unitlength}{3pt}
\psset{unit=3pt}
\psset{runit=2pt}
\psset{linewidth=0.2}
\begin{pspicture}(0,.5)(5,6)
\psline(3,0)(3,2)
\psline(3,2)(1,4)
\psline(3,2)(4,3)
\psline(2,3)(4,5)
\psline(3,4)(2,5)
\end{pspicture}}

\newcommand\troistrois{
\setlength{\unitlength}{3pt}
\psset{unit=3pt}
\psset{runit=2pt}
\psset{linewidth=0.2}
\begin{pspicture}(-0.5,.5)(6.5,6)
\psline(3,0)(3,2)
\psline(3,2)(0.5,4.5)
\psline(1.5,3.5)(2.5,4.5)
\psline(3,2)(5.5,4.5)
\psline(4.5,3.5)(3.5,4.5)
\end{pspicture}}

\newcommand\troisquatre{
\setlength{\unitlength}{3pt}
\psset{unit=3pt}
\psset{runit=2pt}
\psset{linewidth=0.2}
\begin{pspicture}(1,.5)(6,6)
\psline(3,0)(3,2)
\psline(3,2)(5,4)
\psline(3,2)(2,3)
\psline(4,3)(2,5)
\psline(3,4)(4,5)
\end{pspicture}}

\newcommand\troiscinq{
\setlength{\unitlength}{3pt}
\psset{unit=3pt}
\psset{runit=2pt}
\psset{linewidth=0.2}
\begin{pspicture}(1,.5)(7,6)
\psline(3,0)(3,2)
\psline(3,2)(6,5)
\psline(3,2)(2,3)
\psline(4,3)(3,4)
\psline(5,4)(4,5)
\end{pspicture}}

\begin{document}
 
\title{QED Hopf algebras on planar binary trees}
 
\author{Christian Brouder \\
Laboratoire de Min\'eralogie-Cristallographie, CNRS UMR 7590,
 Universit\'es Paris 6 et 7, \\IPGP, 4 place Jussieu,
  75252 Paris Cedex 05,
  France. {\small \tt brouder@lmcp.jussieu.fr}
\and
Alessandra Frabetti \\
Institut Girard Desargues, CNRS UMR 5028, Universit\'e Lyon 1,\\
21 avenue Claude Bernard, 69622 Villeurbanne, France.\\
{\small \tt frabetti@desargues.univ-lyon1.fr}}
\date{\today}
 
\maketitle


\begin{abstract}
In this paper we describe the Hopf algebras on planar binary trees 
used to renormalize the Feynman propagators of quantum 
electrodynamics, and the coaction which describes the renormalization 
procedure. Both structures are related to some semi-direct coproduct 
of Hopf algebras. 
\end{abstract}




\section{Introduction}

Planar binary trees exhibit surprinsigly rich structures, often
related to analogue ones on (non binary) rooted trees. 
For instance, in the last decades several Hopf algebras on families 
of trees have been discovered in different frameworks: by R.Grossman
and R.G.~Larson \cite{GrossmanLarson} in connection to the Butcher
group introduced by J.C. Butcher \cite{Butcher} to solve differential 
equations; by D.~Kreimer \cite{Kreimer98} to describe the
renormalization of perturbative quantum field theory; by J.-L.~Loday 
and M.~Ronco \cite{LodayRonco1} in the framework of dendriform
algebras; by L.~Foissy \cite{Foissy} as a noncommutative extension of 
the Butcher-Kreimer algebra on rooted trees. 
In particular, the relationship between some of these Hopf algebras 
have been studied by F.~Painate \cite{Panaite}, by Foissy
\cite{Foissy} and by R.~Holtkamp \cite{Holtkamp}. 
Moreover, these Hopf algebras present a universal property which was 
first described by A.~Connes and Kreimer in \cite{ConnesKreimerIHES}, 
and then used by I.Moerdijk in \cite{Moerdijk} to introduce a large 
class of Hopf algebras. 

Planar binary trees were also used in \cite{Brouder} to solve 
perturbatively the system of functional differential equations satisfied 
by the 2-points correlation functions of quantum electrodynamics (QED). 
QED is the quantum field theory which describes the dynamics of 
interacting electrons and photons. 
The interaction between these two particles is usually represented by 
Feynman diagrams. In this context, each tree corresponds to a finite sum 
of appropriate Feynman diagrams and the explicit relations are given in 
\cite{BrouderFrabetti1}. 
Following the classical Feynman rules, cf. for instance \cite{ItzyksonZuber} 
or \cite{PeskinSchroeder}, we consider two Feynman amplitudes associated 
to each tree, one for the electron propagator and one for the photon 
propagator. 
These amplitudes are in general sums of divergent integrals which need 
to be renormalized, and in \cite{BrouderFrabetti1} we give a 
perturbative solution in term of trees of the equations 
satisfied by the renormalized propagators of QED. 

In principle, the renormalization procedure can be described directly on 
the propagators as the action of a group, called the renormalization group.  
In practice, the elements of this group are known only through computations 
made on Feynman diagrams: the so-called Forest Formula \cite{ItzyksonZuber} 
describes the relation between the perturbative coefficients of the 
propagators before and after the renormalization. 
In \cite{Kreimer98,ConnesKreimerI}, D.~Kreimer and A.~Connes discovered 
that the operations involved in this formula define the structure of a 
commutative Hopf algebra on the set of Feynman diagrams labeled by 
some indices. 
This result shows that the labeled Feynman diagrams are the natural local 
coordinates of quantum field theory, and that the renormalization group 
can be recovered as the set of characters of this Hopf algebra. 
Therefore, the renormalization procedure is known if we can construct 
the characters from the data we know of the quantum field 
theory: the Feynman amplitudes and the counterterm maps imposed by the 
physical renormalization prescription. 

In scalar field theory, the amplitudes and the counterterms are scalar 
maps which preserve the junction of Feynman graphs\footnote{
We consider Feynman graphs with one external leg removed.
The junction is then the same as the disjoint union.}, 
that is, they are characters of the Connes-Kreimer Hopf algebra 
given in \cite{ConnesKreimerI}. 
The relationship between the renormalization group and the local coordinates 
Hopf algebra is then the classical Tannaka-Krein duality which holds 
between each affine algebraic group and its coordinate ring. 

In vector or spinor valued field theories (such as QED), the 
propagators\footnote{
We consider here the usual propagators multiplied by the inverse 
of the free propagator.} 
are matrices, hence the Feynman amplitudes and the counterterms are maps 
which take value in a non-commutative ring. Because of Feynman rules, 
they still respect the product between Feynman graphs, but they can not be 
usual characters of a commutative coordinate ring. 
In this case, the Connes-Kreimer commutative Hopf algebra dual to the 
renormalization group of course still exists, but it does not help us 
to recover the group through the physical data, because neither the 
amplitudes nor their matrix elements are characters of this algebra. 
The alternative approach is to look for a suitable algebra whose 
matrix-valued characters are the known Feynman amplitudes and 
counterterm maps. It might not exist, or not be a Hopf algebra. 
In fact, these ``non-commutative characters'' do not satisfy any known 
duality principle, for two reasons. 
Given a group $G = \Hom_{Alg}(\H,\A)$ through a set of algebra homomorphisms 
between two non-commutative algebras, the group law on $G$ induces 
a coproduct on $\H$ which, in general, takes value in the free product 
$\H \star \H$, therefore $\H$ is not necessarily a Hopf algebra in the usual 
sense. Dually, if $\H$ is a non-commutative Hopf algebra, then the set 
$\Hom_{Alg}(\H,A)$ is a groupoid, in general, but not necessarily a group. 

In this paper, we show that there exists a non-commutative Hopf algebra 
which represents the renormalization group of QED, and 
an associated coaction on the algebra dual to the propagators which 
describes the renormalization in local coordinates. 
In our case, we expand the QED propagators as asymptotic series over the 
planar binary trees instead of the Feynman diagrams, therefore our local 
non-commutative coordinates are the trees. 
The coproducts and coactions are forced by the relationship between the bare 
and the renormalized propagators found in \cite{BrouderFrabetti1}. 
The resulting Hopf algebras of renormalization for photons and electrons 
look very different, both for the algebra and the coalgebra structures. 
However, they can be interpretated as semidirect coproducts of similar 
Hopf algebras, and thus directly related to a standard form 
of the renormalization group. 
\bigskip

The paper is organized as follows. In the second section we recall the 
algebraic tools needed to present a non-commutative version of the 
the renormalization group, which is a semidirect product of two groups, 
and of the renormalization action. 
The main tools come from the semidirect or smash coproduct of Hopf algebras, 
introduced by R.~Molnar in \cite{Molnar}.  

In the third, fourth and fifth sections, we define the non-commutative 
Hopf algebras which correspond to the electron and photon propagators; 
the Hopf algebra which corresponds to the renormalization of the coupling 
constant of QED; and finally the renormalization Hopf algebras and the 
renormalization coactions for the electron and for the photon propagators. 
To describe these structures we only need some grafting and pruning 
operations on trees. The choice of such operations, which looks 
apparently arbitrary, is in fact forced by the combinatorial operations 
on the Feynman graphs related to the trees, cf. \cite{BrouderFrabetti1}.  
It is then even more surprising how the basic operations on trees 
turn out to be deeply related to those used by J.-L.~Loday in his 
{\em Arithmetree}, \cite{LodayATree}. 

The main application of these Hopf algebras, namely the renormalization 
of QED propagators, is recalled in the last section. 
\bigskip

\noindent
{\bf Notations.} 
We suppose that all vector spaces and algebras are defined over the  
field $\C$ of complex numbers, but this choice is not necessary. 
For any set $X$, we denote by $\C X$ the vector space spanned by $X$, 
by $\C\langle X \rangle$ the tensor algebra on $X$ (noncommutative 
polynomials), and by $\C[X]$ the symmetric algebra on $X$ (commutative 
polynomials). 
\bigskip

\noindent
{\bf Aknowledgments.} 
A.~Frabetti warmly thanks the Swiss National Foundation for Scientific 
Research for the support while the paper was written, and the Mathematics 
Departement in Lausanne University for their warm hospitality. 
This is IGPG contribution \#0000.


\section{Renormalization group and semidirect coproduct of Hopf algebras} 
\label{semidirect}

The renormalization of quantum fields can be formalized as an action 
of the renormalization group on the set of propagators associated 
to the quantum fields: a bare propagator $D$ is turned into a renormalized 
propagator $\bar D = D \cdot \eta$ by the action of an element $\eta$ 
of the renormalization group. 
In perturbative quantum field theory, all these groups and sets are 
made of formal series in the powers of the coupling constant, 
which is the fine structure constant $\alpha$ in QED (i.e. the square 
of the electric charge divided by $4\pi$). 
Among such series, two basic operations are possibly allowed and 
determine the algebraic part of the renormalization procedure: 
the (pointwise) product and the composition, or substitution. 

The set of propagators is in fact a group $G^p$ of series with 
the pointwise product, since their constant term is invertible. 
The composition, instead, is the natural operation 
in the group $G^c$ which renormalizes the coupling constant. 
Such series have zero constant term, and invertible first term. 
By substitution, the group $G^c$ also acts on $G^p$ from the right, 
and the action $(f,\varphi) \mapsto f^\varphi$ is associative with 
respects to the composition in $G^c$, i.e.
\begin{eqnarray*}
	f^{(\varphi \psi)} &=& (f^{\varphi})^{\psi}, 
\end{eqnarray*} 
and commutes with the product in $G^p$, in the sense that 
\begin{eqnarray*}
	(f g)^\varphi &=& (f^\varphi) (g^\varphi).  
\end{eqnarray*} 

Then, in QED we can distinguish two renormalization groups, one for the 
electron and one for the photon propagators.  

\begin{ale}{\bf The QED renormalization groups.}
\label{renormalization group} 
The electron renormalization group is the semidirect product 
$G^c \ltimes G^p$ made of pairs $(\varphi,f)$ in the direct product 
$G^c \times G^p$ with group law 
\begin{eqnarray*}
	(\varphi,f) \ldot (\psi,g) &:=& 
	(\varphi \psi, f^{\psi} g). 
\end{eqnarray*} 
The renormalization procedure is the right action of $G^c \ltimes G^p$ 
on $G^p$ obtained by embedding $G^p$ in $G^c \ltimes G^p$ through 
$f \mapsto (1_c,f)$, applying the semidirect product law in 
$G^c \ltimes G^p$ and then projecting onto the $G^p$ component, that is, 
\begin{eqnarray*}
	f \ldot (\varphi,f) &:=& f^{\varphi} f.
\end{eqnarray*}

In this model, the $G^c$ component of $G^c\ltimes G^p$ represents 
the renormalization of the fine structure constant, while the $G^p$ 
component represents the inverse $Z_2^{-1}$ of the electron 
renormalization factor. 
In the analogue situation for the photon renormalization, it was proved 
by J.C.~Ward \cite{Ward} that the fine structure constant is renormalized 
exactly by the inverse $Z_3^{-1}$ of the photon renormalization factor. 
In other words, we identify the series of type $G^c$ and $G^p$ by 
multiplying or dividing by the fine structure constant. 
In our model, this identification consists of a map 
$s : G^c \longrightarrow G^p$ which is a 1-cocycle of $G^c$ with values 
in $G^p$, that is,  
\begin{eqnarray*}
	s(\psi) [s(\varphi \psi)]^{-1} [s(\varphi) \psi] 
	&=& 1_p \qquad\mbox{for all $\varphi,\psi \in G^c$}. 
\end{eqnarray*}

The photon renormalization group is then the group $G^c$ itself, 
and the renormalization procedure simply becomes the action on $G^p$ 
of $G^c$ identified with the subgroup $G^c \ltimes s(G^c)$ of 
$G^c \ltimes G^p$, 
\begin{eqnarray*}
	f \sdot \varphi &:=& f^{\varphi} s(\varphi). 
\end{eqnarray*}
\end{ale}

\begin{ale}{\bf Feynman amplitudes and characters.}
\label{characters} 
The Feynman bare and renormalized amplitudes $U$, $R$, and the counterterm 
maps $C$ are all the data which allow to reconstruct the bare and the 
renormalized propagators and the elements of the renormalization group,  
starting from the appropriate set of Feynman diagrams (or trees in our case). 

Moreover, by definition of the Feynman rules, they preserve the natural 
product which joins together two (amputated) Feynman diagrams. 

For scalar field theories, these maps take scalar values, so they 
can be recognized as characters of the coordinate rings 
of the groups involved. 
In fact, since $G^c$ and $G^p$ are two affine groups, the semidirect product 
$G^c \ltimes G^p$ is also affine. Let us denote by $\C(G^c)$, $\C(G^p)$ and 
$\C(G^c \ltimes G^p)$ their coordinate rings. They are commutative Hopf 
algebras, in perfect duality with the original groups. 
More precisely, the groups can be reconstructed as the sets 
\begin{eqnarray*}
	&& G^c \cong \Hom_{Alg}(\C(G^c),\C), 
	\quad G^p \cong \Hom_{Alg}(\C(G^p),\C), 
	\quad G^c \ltimes G^p \cong \Hom_{Alg}(\C(G^c \ltimes G^p),\C), 
\end{eqnarray*}
of characters of the coordinate rings, which are the algebra homomorphisms 
from the rings to the field of scalars, endowed with the convolution 
products. 

For QED, the maps $U$, $R$ and $C$ do not anymore have scalar values.
For single particle Green functions, the maps $U$ and $R$ take value in the 
ring of $4 \times 4$ complex matrices, and $C$ is defined by the value of 
$U$ and $R$ at some fixed momentum, through the renormalization conditions. 
Therefore, $C$ is {\em a priori} a matrix. 
In the case of QED in flat space time, Lorentz invariance implies that $C$ 
is a scalar multiplied by a fixed matrix. However, for applications in 
curved space time, or to treat several fermions at once, it is interesting 
to allow matrix-valued counterterms. 
Therefore, in general, the maps $U$, $R$ and $C$ are not anymore characters 
of the coordinate rings of the groups. A duality, if it exists, should then 
be searched between the groups $G^p$, $G^c$, $G^c \ltimes G^p$ and some 
algebras $\H(G^p)$, $\H(G^c)$, $\H(G^c \ltimes G^p)$ such that 
\begin{eqnarray*}
	&& U, R \in \Hom_{Alg}(\H(G^p),\A), 
	\qquad C \in \Hom_{Alg}(\H(G^c \ltimes G^p),\A), 
\end{eqnarray*}
where $\A$ is the non-commutative ring where $U$, $R$ and $C$ 
take values. 
This leads us to consider some non-commutative versions of the 
coordinate rings. 
\end{ale}

\begin{ale}{\bf The Hopf algebra of a semidirect product of groups.}
\label{semidirect group} 
Let us recall how to construct the coordinate ring $\C(G^c \ltimes G^p)$ 
and the coactions on $\C(G^p)$. 
Denote by $\D^c : \C(G^c) \longrightarrow \C(G^c) \otimes \C(G^c)$ and 
$\D^p : \C(G^p) \longrightarrow \C(G^p) \otimes \C(G^p)$ the coproducts 
dual to the group laws of $G^c$ and $G^p$, in the sense that 
\begin{eqnarray*}
	\< \varphi \psi, a \> = \< \varphi \otimes \psi, \D^c a \> 
	&\mbox{and}& 
	\< f g, b \> = \< f \otimes g, \D^p b \>,
\end{eqnarray*}  
where $\<\ ,\ \>: G^c \times \C(G^c) \longrightarrow \C$ is the 
evaluation map $\< \varphi, a \>= a(\varphi)$. 
Also, denote by $\dd : \C(G^p) \longrightarrow \C(G^p) \otimes \C(G^c)$ 
the coaction dual to the action of $G^c$ on $G^p$, 
\begin{eqnarray*}
	\< h^\varphi, b \> &=& \< h \otimes \varphi, \dd (b) \> . 
\end{eqnarray*}  
The map $\dd$ is coassociative with respect to $\D^c$ and commutes 
with $\D^p$. 

Then, the coordinate ring $\C(G^c \ltimes G^p)$ is the tensor product 
algebra $\C(G^c) \otimes \C(G^p)$, endowed with the coproduct $\Dl$ dual 
to the group law $\ldot$, i.e. 
\begin{eqnarray*}
	\< (\varphi,f) \ldot (\psi,g), a \otimes b \> &=& 
	\< \varphi \otimes f \otimes \psi \otimes g, \Dl (a \otimes b) \>. 
\end{eqnarray*}  
Explicitly, $\Dl$ is the algebra morphism given by 
\begin{eqnarray*}
	 \Dl (a \otimes b) &=& 
	\D^c(a) \ [(\dd \otimes \Id) \D^p(b)], 
\end{eqnarray*}
where we omit the symbol of the componentwise product in the algebra 
$\C(G^c) \otimes \C(G^p)$ between the image of $\D^c$ in 
$\C(G^c) \otimes 1 \otimes \C(G^c) \otimes 1$, and the image of 
$(\dd \otimes \Id) \D^p$ in 
$1 \otimes \C(G^p) \otimes \C(G^c) \otimes \C(G^p)$. 

Moreover, the action of $G^c \ltimes G^p$ on $G^p$ induces a dual coaction 
$\ddl: \C(G^p) \longrightarrow \C(G^p) \otimes \C(G^c \ltimes G^p)$, 
which is simply the second component of the coproduct, i.e. 
\begin{eqnarray*}
	 \ddl (b) &=& (\dd \otimes \Id) \D^p(b).  
\end{eqnarray*}

Similarly, if we denote by $\sigma : \C(G^p) \longrightarrow \C(G^c)$ the 
linear map dual to the 1-cocycle $s:G^c \longrightarrow G^p$, 
then $\sigma$ satisfies the identity 
\begin{eqnarray*}
	(m^c_{24} \otimes m^c_{135}) 
	(\Id \otimes \D^c \otimes \Id \otimes \Id) 
	(\sigma \otimes \sigma \otimes \sigma \otimes \Id) 
	(\Id \otimes S^p \otimes \dd) (\D^p)^2 
	&=& \D^c \ i^c \ \epsilon^p , 
\end{eqnarray*} 
where $m^c_{ijk}$ denotes the multiplication in $\C(G^c)$ applied to positions 
$(i,j,k)$, $S^p$ is the antipode of the Hopf algebra $\C(G^p)$, 
$(\D^p)^2 = (\D^p \otimes \Id) \D^p = (\Id \otimes \D^p) \D^p$, 
$i^c : \C \longrightarrow \C(G^c)$ is the unit of $\C(G^c)$, 
and finally $\epsilon^p : \C(G^p) \longrightarrow \C$ is the counit 
of $\C(G^p)$. This condition is equivalent to require that 
\begin{eqnarray*}
	\D^c \sigma &=& (\sigma \otimes \Id) 
	(\Id \otimes m^c) (\dd \otimes \sigma) \D^p. 
\end{eqnarray*}

The coaction of $\C(G^c)$ on $\C(G^p)$ dual to the action $\sdot$ is then 
the algebra morphism $\dds: \C(G^p) \longrightarrow \C(G^p) \otimes \C(G^c)$ 
given by 
\begin{eqnarray*}
	\dds(b) &=& (\Id \otimes m^c) (\dd \otimes \sigma) \D^p (b). 
\end{eqnarray*}
\end{ale}

\begin{ale}{\bf The semidirect coproduct of Hopf algebras.} 
\label{semidirect Hopf}
The formulas of section (\ref{semidirect group}) make sense for all Hopf 
algebras, even not commutative ones. The generalisation to arbitrary Hopf 
algebras has been studied by R.~Molnar \cite{Molnar}, B.~Lin \cite{Lin}, 
D.~Radford \cite{Radford}, S.~Majid \cite{Majid} and others. 

Let $\H^c$ and $\H^p$ be two Hopf algebras with multiplications $m^c$, $m^p$ 
and coproducts $\D^c$, $\D^p$. Suppose that $\H^c$ coacts on $\H^p$ from 
the right, and that the coaction $\dd: \H^p \longrightarrow \H^p \otimes \H^c$ 
satisfies 
\begin{eqnarray}
	(\dd \otimes \Id) \dd &=& (\Id \otimes \D^c) \dd 
	\label{D1-coassociative} \\ 
	(\D^p \otimes \Id) \dd &=& m_{24}^3 (\dd \otimes \dd) \D^p , 
	\label{D2-commute} 
\end{eqnarray}
where $m_{24}^{3}$ multiplies what is in the position $2$ by what is in 
the position $4$ and puts it in the position $3$. 
Then, the semidirect or smash coproduct $\H^c \ltimes \H^p$ is the 
tensor algebra $\H^c \otimes \H^p$ endowed with the coproduct 
\begin{eqnarray}
\label{semidirect coproduct}
	 \Dl (a \otimes b) &:=& \D^c(a) \ 
	[(\dd \otimes \Id) \D^p(b)] , 
	\qquad a \in \H^c, b \in \H^p, 
\end{eqnarray}
and the counit 
$\epsilon^{\ltimes}(a \otimes b) := \epsilon_1(a) \epsilon_2(b)$. 

Molnar proved in \cite{Molnar} that $\H^c \ltimes \H^p$ is a coalgebra. 
In particular, it follows that the map 
$\ddl: \H^p \longrightarrow \H^p \otimes (\H^c \ltimes \H^p)$
given by $\ddl(b)= (\dd \otimes \Id) \D^p(b)$ is a coaction, 
i.e. it is coassociative with respect to $\Dl$. 

He also proved that $\H^c \ltimes \H^p$ is a bialgebra if $\H^c$ is 
commutative. In this case it is also a Hopf algebra, with antipode 
\begin{eqnarray*}
	S^{\ltimes}(a \otimes b) &:=& S^c(a) \ 
	[\tau (\Id \otimes S^c) \dd (S^p b)] 
	= \tau (S^p \otimes S^c) (\Id \otimes m^c) (\dd \otimes \Id)
	(b \otimes a). 
\end{eqnarray*} 
Moreover, in this case the coaction $\ddl$ is also an algebra morphism.
\end{ale}

\begin{ela}{\bf Lemma.} 
\label{sigma}
Let $\H^c$ and $\H^p$ be two Hopf algebras such that $\H^c$ coacts on 
$\H^p$ as above ($\H^c$ is not necessarily commutative). 
Suppose that there exists a map $\sigma : \H^p \longrightarrow \H^c$ 
with the property that if $\dds: \H^p \longrightarrow \H^p \otimes \H^c$ 
is the map defined by  
\begin{eqnarray*}
	\dds &:=& m_{23}^2 \ (\dd \otimes \sigma) \D^p 
\end{eqnarray*}
then $\sigma$ interwines $\dds$ and $\D^c$, i.e. 
\begin{eqnarray}
\label{interwine}
	\D^c \sigma &=& (\sigma \otimes \Id) \dds. 
\end{eqnarray} 
Then $\dds$ is coassociative with respect to $\D^c$. 
\end{ela}

\begin{proof}
Let us adopt the following Sweedler conventions:
\begin{eqnarray*}
	&& \D^c(a) = \sum a_{(1)} \otimes a_{(2)}, \qquad 
	\D^p(b) = \sum b_{(1)} \otimes b_{(2)}, \qquad 
	\dd(b) = \sum b_{(l)} \otimes b_{(r)}.  
\end{eqnarray*}
Then for any $b \in \H^p$ we have 
\begin{eqnarray*}
	(\dds \otimes \Id) \dds (b) &=& \sum
	\dds (b_{(1l)}) \otimes b_{(1r)} \sigma (b_{(2)}) \\ 
	&=& \sum b_{(1l1l)} \otimes b_{(1l1r)} \sigma (b_{(1l2)}) 
	\otimes b_{(1r)} \sigma (b_{(2)}) \\
	&\stackrel{(1)}{=}& \sum
	b_{(11ll)} \otimes b_{(11lr)} \sigma (b_{(12l)}) 
	\otimes b_{(11r)} b_{(12r)} \sigma (b_{(2)}) \\
	&\stackrel{(2)}{=}& \sum
	b_{(11l)} \otimes b_{(11r1)} \sigma (b_{(12l)}) 
	\otimes b_{(11r2)} b_{(12r)} \sigma (b_{(2)}) \\
	&\stackrel{(3)}{=}& \sum
	b_{(1l)} \otimes b_{(1r1)} \sigma (b_{(21l)}) 
	\otimes b_{(1r2)} b_{(21r)} \sigma (b_{(22)}) \\
	&\stackrel{(4)}{=}& \sum
	b_{(1l)} \otimes b_{(1r1)} \left(\sigma (b_{(2)})\right)_{(1)} 
	\otimes b_{(1r2)} \left(\sigma (b_{(2)})\right)_{(2)} \\
	&=& \sum b_{(1l)} \otimes \D^c \left( b_{(1r)} \sigma(b_{(2)}) \right) 
	= (\Id \otimes \D^c) \dds (b), 
\end{eqnarray*}
where the equality $(1)$ follows from (\ref{D2-commute}) applied 
to $b_{(1)}$, the equality $(2)$ follows from (\ref{D1-coassociative}) 
applied to $b_{(11)}$, the equality $(3)$ follows from the coassociativity 
of $\D^p$ applied to $b$, and the equality $(4)$ follows from 
(\ref{interwine}) applied to $b_{(2)}$. 
\end{proof}


\section{Propagators Hopf algebras on trees}
\label{propagators}

\begin{ale}{\bf Planar binary trees.} 
By {\em planar binary tree} we mean a connected and oriented planar graph 
with no cycle, such that each internal vertex has one incoming and two 
outgoing edges. The incoming and outgoing edges of a tree are called 
respectively the {\em root} and the {\em leaves}. 
Such trees are naturally graded by the number of internal vertices, 
that we call the {\em order}. We denote by $|t|$ the order of a tree $t$. 
Up to continuous transformations of the plane which fix the root and 
the leaves, there are $c_n = \frac{(2n)!}{n!(n+1)!}$ trees with order $n$. 
We denote by $Y_n$ the set of trees $t$ with $|t|=n$, and by 
$Y = \bigcup_{n \geq 0} Y_n$ the set of all planar binary trees. 
Here are the sets of trees with order $0$, $1$, $2$ and $3$: 
\begin{eqnarray*} 
	Y_0 &=& \{ \| \}, \\ 
	Y_1 &=& \{ \Y \}, \\
	Y_2 &=& \{ \deuxun, \deuxdeux \}, \\ 
	Y_3 &=& \{ \troisun, \troisdeux, \troistrois, \troisquatre, 
		\troiscinq \}. 
\end{eqnarray*}

Let $\vee: Y_n \times Y_m \longrightarrow Y_{n+m+1}$ denote the 
map which grafts two trees on a new root, for instance, 
\begin{eqnarray*}
	\Y \vee \Y = \troistrois , &\qquad& 
	\deuxdeux \vee \| = \troisdeux . 
\end{eqnarray*}
Then, each tree $t \neq \|$ is the grafting $t = t^l \vee t^r$ 
of two uniquely determined trees $t^l, t^r$ with smaller order. 
\end{ale}

\begin{ale}{\bf The products {\em over} and {\em under}.} 
\label{pro} 
Following the notations of J.-L.~Loday and M.~Ronco in \cite{LodayATree}, 
\cite{LodayRonco2}, we call {\em over} and {\em under} the graded 
products $\lpro, \rpro: Y_n \times Y_m \longrightarrow Y_{n+m}$ 
defined by the recurrence relations
\begin{eqnarray*}
	t \lpro s &:=& (t \lpro s^l) \vee s^r 
	\quad\mathrm{for}\quad s=s^l\vee s^r, \\
	t \lpro \| &:=& t, 
\end{eqnarray*}
and similarly 
\begin{eqnarray*}
	t \rpro s &:=& t^l \vee (t^r \rpro s) 
	\quad\mathrm{for}\quad t=t^l\vee t^r, \\
	\| \rpro s &:=& s. 
\end{eqnarray*}
These operations graft one tree on the other one 
according to the rules $t \lpro s = \lgraft{s}{t}$ and 
$t \rpro s = \rgraft{t}{s}$. 
For instance, 
\begin{eqnarray*}
	\deuxdeux \lpro \Y = \troisdeux , &\qquad& 
	\Y \lpro \deuxdeux = \troistrois, \\ 
	\deuxun \rpro \Y = \troistrois , &\qquad& 
	\Y \rpro \deuxun = \troisquatre . 
\end{eqnarray*}
Both products are clearly associative (non commutative), and for both
the root tree $\|$ is a unit. Moreover, any tree 
$t=t^l \vee t^r$ can be decomposed as $t=t^l \lpro (\|\vee t^r)$ 
or as $t = (t^l \vee \|) \rpro t^r$. 
Hence, the trees of the form $\| \vee t =: V(t)$, for any $t\in Y$, 
form a system of generators of $(Y,\lpro)$, and similarly 
the trees of the form $t \vee \|$ form a system of generators 
of $(Y,\rpro)$. 
\end{ale}

\begin{ale}{\bf The pruning coalgebras.}
\label{pruning coalgebras}
Identify $\C Y$ with its linear dual $\C Y^*$, and consider 
the coproducts $\DPp, \DPe :\C Y \longrightarrow \C Y \otimes \C Y$ 
dual of the products $\lpro$ and $\rpro$ respectively, 
\begin{eqnarray*}
	\DPp(t) &=& \sum_{t=t_1 \lpro t_2} t_1 \otimes t_2, \\ 
	\DPe(t) &=& \sum_{t=t_1 \rpro t_2} t_1 \otimes t_2. 
\end{eqnarray*}
Of course, $\DPp$ and $\DPe$ are graded coassociative operations, 
and together with the counit $\epsilon$ dual to the unit $\|$, 
defined as $\epsilon(\|)=1$ and $\epsilon(t)=0$ if $t \neq \|$, 
they define on $\C Y$ two different structures of graded coalgebra. 

The coproducts $\DPp$ and $\DPe$ break all the branches of a tree which are 
respectively on the left and on the right of the root, and places them 
on the same side. 
It is useful to give a recursive definition of these coporducts: 
for any $t,s \in Y$ we have 
\begin{eqnarray}
	\DPp(\|) &=& \| \otimes \|, \nonumber \\ 
\label{left pruning coproduct}
	\DPp(t \vee s) &=& t \vee s \otimes \| 
	+ \sum_{\DPp t} t_{(1)} \otimes t_{(2)} \vee s,  
\end{eqnarray}
and similarly 
\begin{eqnarray}
	\DPe(\|) &=& \| \otimes \|, \nonumber \\ 
\label{right pruning coproduct}
	\DPe(t \vee s) &=& \| \otimes t \vee s 
	+ \sum_{\DPe s} t \vee s_{(1)} \otimes s_{(2)},  
\end{eqnarray}
where we use the standard Sweedler notation 
$\DPp(t)=\sum t_{(1)} \otimes t_{(2)}$ and 
$\DPe(s)=\sum s_{(1)} \otimes s_{(2)}$. 
The pruning operator of \cite{Brouder} is the reduced coproduct 
$P(t) = \DPe(t) - t \otimes \| - \| \otimes t$. 
\end{ale}

\begin{ale}{\bf The photon and electron propagator Hopf algebras.} 
\label{pruning Hopf}
If we extend the pruning coproducts $\DPp$ and $\DPe$ multiplicatively 
on tensor products of trees, and we set the root tree $\|$ as unit, 
we obtain two different Hopf algebras $\Hp$ and $\He$, which are 
neither commutative nor cocommutative. Therefore we set 
$\Hp,\He :=\C\langle Y \rangle /(1- \| )$ 
as the free associative algebras on the set of trees 
where we identify the formal unit $1$ with the root tree $\|$, 
and we consider $\Hp$ with the Hopf structure induced by $\DPp$, 
and $\He$ with the Hopf structure induced by $\DPe$. 
For notational convenience, we omit the tensor product symbols. 

Beside the natural grading coming from the tensor powers, 
on a tensor product of trees we can define 
a total order as the sum of the orders of the trees, 
\begin{eqnarray*}
	|t_1 \ldots t_k| &=& |t_1| +\cdots+ |t_k|. 
\end{eqnarray*}
Then the algebras $\Hp$ and $\He$ are graded connected Hopf algebras, 
with homogeneous components 
\begin{eqnarray*}
	\Hp_n, \He_n &=& \bigoplus_{n_1+\cdots+n_k=n} 
	\C Y_{n_1} \otimes\dots\otimes \C Y_{n_k}. 
\end{eqnarray*}

In particular, the electron pruning antipode $\SPe$ is the graded algebra 
anti-morphism automatically defined on generators by the recursive 
formula $\SPe(\|)=\|$ and
\begin{eqnarray*}
	\SPe(t) &=& -t -\sum_{P(t)}  \SPe(t_{(1)}) t_{(2)} 
	= -t -\sum_{P(t)}  t_{(1)} \SPe(t_{(2)}).  
\end{eqnarray*}
Since $\SPe$ plays an explicit role in the renormalization of the 
electron propagator, we give a few examples: 
\begin{eqnarray*}
	\SPe(\Y) &=& -\Y , \\ 
	\SPe(\deuxdeux) &=& -\deuxdeux+\Y^2, \\ 
	\SPe(\deuxun) &=& -\deuxun, \\ 
	\SPe(\troiscinq) &=& -\troiscinq + \deuxdeux \Y
	+ \Y \deuxdeux - \Y^3, \\ 
	\SPe(\troisquatre) &=& - \troisquatre + \Y \deuxun . 
\end{eqnarray*}
Notice that the coproduct $\DPe$ is neither commutative 
nor cocommutative, and $\SPe \circ \SPe \neq \Id$.
\end{ale}


\section{Charge Hopf algebra on trees}
\label{charge} 

\begin{ale}{\bf The charge algebra.} 
\label{algebraHa}
Let $\Ha := \C[V(t), t \in Y]$ be the polynomal algebra generated 
by all trees of the form $V(t)=\| \vee t$. 
Since each tree $t \in Y$ can be uniquely decomposed as 
$t = t_l \lpro V(t_r)$, the map $V(t) \mapsto V(t)$ and $1 \mapsto \|$ 
is an algebra isomorphism from $\Ha$ to the abelianization of $(\C Y,\lpro)$. 
Under the inverse of this isomorphism, the natural homogeneous component 
$\C Y_n$ of degree $n$ in $\C Y$ corresponds to the subspace 
$\Ha_n = \bigoplus_{n_1\leq \cdots \leq n_k} 
\C V(Y_{n_1}) \otimes\cdots\otimes \C V(Y_{n_k})$ 
of total degree $n=n_1+...+n_k+k$ in $\Ha$. 

From now on, we identify $\Ha$ with $(\C Y,\lpro)_{ab}$, and represent 
the unit $1$ as the root tree $\|$. 
\end{ale}

\begin{ale}{\bf The charge Hopf algebra.} 
\label{HopfHa}
Define a coproduct $\Da: \Ha \longrightarrow \Ha \otimes \Ha$ and 
a coaction $\dd : \Ha \longrightarrow \Ha \otimes \Ha$ as the two 
linear operators satisfying the following recursive relations:  
\begin{eqnarray*}
	\Da \| &=& \| \otimes \|, \\ 
	\Da V(t) &=& \| \otimes V(t) + \dd V(t), \\ 
	\Da (t \vee s) &=& \Da t \lpro \Da V(s) ; 
\end{eqnarray*}
and 
\begin{eqnarray*}
	\dd \| &=& \| \otimes \|, \\  
	\dd V(t) &=& (V \otimes \Id) \dd (t), \\  
	\dd (t \vee s) &=& \Da t \lpro \dd (V(s)) .  
\end{eqnarray*}
For instance, the coproduct on small generator trees yields 
\begin{eqnarray*}
	\Da \Y &=& \Y\otimes \| + \| \otimes \Y , \\ 
	\Da \deuxdeux &=& \deuxdeux\otimes \| + \| \otimes \deuxdeux , \\ 
	\Da\troisquatre &=& \troisquatre \otimes \| 
	+ \deuxdeux \otimes \Y + \| \otimes \troisquatre ,\\
	\Da\troiscinq &=& \troiscinq \otimes \| + \| \otimes \troiscinq .
\end{eqnarray*}
Similarly, the coaction on small generator trees yields 
\begin{eqnarray*}
	\dd \Y &=& \Y \otimes \| , \\ 
	\dd \deuxdeux &=& \deuxdeux \otimes \| , \\ 
	\dd\troisquatre &=& \troisquatre \otimes \| +\deuxdeux \otimes \Y,\\
	\dd\troiscinq &=& \troiscinq \otimes \|.
\end{eqnarray*}

Let $\epsilon : \Ha \longrightarrow \C$ be the linear map 
which sends all the trees to $0$ except the root 
$\|$ which is sent to $1$. 
\end{ale}

\begin{ela}{\bf Theorem.} 
\label{theoremHa}
The algebra $\Ha$ is a graded connected commutative Hopf algebra. 
Moreover, $\dd$ is a right $\Da$-coaction, that is 
\begin{eqnarray*}
	(\dd \otimes \Id ) \dd  &=& (\Id \otimes \Da) \dd . 
\end{eqnarray*}
\end{ela}
 
\begin{proof} 
We first observe that the coproduct preserves 
the grading of $\Ha$, that is 
\begin{eqnarray*}
	\Da(\Ha_n) \subset \bigoplus_{p+q=n} \Ha_p \otimes \Ha_q.  
\end{eqnarray*}
Since $\Ha_0$ is spanned by a single tree $\|$, the graded algebra 
$\Ha$ is connected. 

By recursion arguments, it is then easy to see that the only terms 
of total degree $(n,0)$ and $(0,n)$, in the image of $\Da$, 
consist of the primitive part $t \otimes \|$ and $\| \otimes t$ 
for any tree $t$.
Then, the map $\epsilon$ is a counit for $\Da$, and the antipode 
$\Sp: \Ha \longrightarrow \Ha$ is the graded algebra isomorphism 
automatically defined on the generators by the recursive formula 
\begin{eqnarray}
\label{antipodeHc} 
	\Sp(t) = -t - \sum_{\bar{\Da}(t)} \Sp(t_{(1)}) \lpro t_{(2)},  
\end{eqnarray}
where $\bar{\Da}(t)=\Da(t) - t \otimes \| - \| \otimes t$ is the 
reduced coproduct. 

First we prove by induction that the operator $\dd $ defines a left
$\Da$-coaction of $\Ha$ on itself. 
It is true on $t=\|$. Suppose that it is true 
for all the trees up to order $n$, and let $V(t)$ has order $n+1$. 
Then 
\begin{eqnarray*} 
	(\dd \otimes \Id ) \dd (V(t)) &=& 
	(\dd \circ V \otimes \Id) \dd (t) 
	= (V \otimes \Id \otimes \Id) (\dd \otimes \Id ) \dd (t) \\ 
	&=& (V \otimes \Id \otimes \Id) (\Id \otimes \Da) \dd (t) 
	= (V \otimes \Da) \dd (t) \\ 
	&=& (\Id \otimes \Da) \dd (V(t)) . 
\end{eqnarray*}

Now let $s \vee t = s \lpro V(t)$ has order $n+1$, with $s \neq \|$. 
Then both $s$ and $V(t)$ have order smaller or equal to $n$. 
Let us fix the Sweedler notations 
\begin{eqnarray*}
	\Da(s) = \sum s_{(1)} \otimes s_{(2)}, &\quad&
	\dd(t) = \sum t_{(l)} \otimes t_{(r)}. 
\end{eqnarray*}
On one side we have 
\begin{eqnarray*} 
	(\dd \otimes \Id ) \dd (s \lpro V(t)) &=& 
	(\dd \otimes \Id ) [\Da(s) \lpro \dd V(t) ]  
	=(\dd \otimes \Id) \sum_{\dd t, \Da s} 
	s_{(1)} \lpro V(t_{(l)}) \otimes s_{(2)} \lpro t_{(r)} \\
	&=& \sum_{\dd t, \Da s} 
	\Da(s_{(1)}) \lpro \dd ( V(t_{(l)})) \otimes s_{(2)} \lpro t_{(r)} \\
	&=& [(\Da \otimes \Id)\Da(s)] \lpro [(\dd \otimes \Id) \dd V(t)] , 
\end{eqnarray*}
and on the other side we have 
\begin{eqnarray*} 
	(\Id \otimes \Da) \dd (s \lpro V(t)) &=& 
	[(\Id \otimes \Da) \Da(s)] \lpro [(\Id \otimes \Da) \dd V(t)] ,  
\end{eqnarray*}
so the equality holds by inductive hypothesis. 

Now we prove by induction that the operator $\Da$ is coassociative, 
that is $(\Id \otimes \Da) \Da = (\Da \otimes \Id) \Da$, 
using the fact that $\dd $ is a coaction. Since $\Da$ is multiplicative, 
we only need to prove it on the generators $V(t)$. It is true on $t=\|$. 
Suppose that $\Da$ is coassociative on all the trees with order 
up to $n$, and let $V(t)$ be a generator with order $n+1$. 
Then by definition of $\Da$ we have on one side 
\begin{eqnarray*}
	(\Id \otimes \Da) \Da V(t) &=& 
	\| \otimes \Da V(t) + (\Id \otimes \Da) \dd (V(t)) \\ 
	&=& \| \otimes \| \otimes V(t) + \|\otimes \dd (V(t))
	+ (\Id \otimes \Da) \dd (V(t)), 
\end{eqnarray*}
and on the other side 
\begin{eqnarray*}
	(\Da \otimes \Id) \Da V(t) &=& 
	\Da(\|)\otimes V(t) +  (\Da \otimes \Id) \dd (V(t)) \\  
	&=& \| \otimes \| \otimes V(t)
	+ (\Da \circ V \otimes \Id) \dd (t) \\
	&=& \| \otimes \| \otimes V(t)
	+ (\Id \otimes V \otimes \Id) (\|\otimes\dd (t))
	+ (\dd \circ V \otimes \Id) \dd (t) \\ 
	&=& \| \otimes \| \otimes V(t)
	+ \|\otimes(V \otimes \Id) \dd (t)
	+ (\dd \otimes \Id ) (V \otimes \Id) \dd (t) \\ 
	&=& \| \otimes \| \otimes V(t)
	+ \| \otimes \dd (V(t))
	+ (\dd \otimes \Id ) \dd (V(t)).   
\end{eqnarray*}
Then, the two sides are equal because 
$(\Id \otimes \Da) \dd (V(t)) = (\dd \otimes \Id ) \dd (V(t))$. 
\end{proof}

\begin{ale}{\bf The non-commutative charge Hopf algebra.} 
\label{algebraHaNC} 
Let $\Hat := \C \< V(t), t\in Y \>$ be the algebra of non 
commutative polynomials on the trees of the form $V(t)$. 
Then the charge algebra $\Ha$ is the abelian quotient of $\Hat$. 
Moreover, the isomorphism 
$\Ha \stackrel{\sim}{\longrightarrow} (\C Y,\lpro)_{ab}$ 
of (\ref{algebraHa}) can be lifted to an isomorphism 
$\Hat \stackrel{\sim}{\longrightarrow} (\C Y,\lpro)$. 
Therefore, the formulas employed in (\ref{HopfHa}) to define 
a coproduct $\Da$ and a coaction $\dd$ on $\Ha$ can be adopted 
to define some lifted maps $\Dat$ and $\ddt$ 
from $\Hat$ to $\Hat \otimes \Hat$. 
These lifted maps are defined as the original ones on the generators, 
and no ambiguity comes from a product of generator trees if we 
require $\Dat$ and $\ddt$ to be algebra 
morphisms. 
\end{ale}

\begin{ela}{\bf Theorem.} 
\label{theoremHaNC}
The algebra $\Hat$ is a graded connected Hopf algebra, 
which is neither commutative nor cocommutative. 
\end{ela}
 
\begin{proof} 
We can repeat the proof of (\ref{theoremHa}), since we never used 
the commutativity of the product in $\Ha$. 
\end{proof}


\section{QED Hopf algebra and coactions on trees}

\begin{ale}{\bf The electron and photon coactions.} 
\label{coactiondd}
Since $\Hat \cong \C Y$ as a vector space, the coaction $\ddt$ 
on $\Hat$ given in (\ref{algebraHaNC}) can be seen as a linear map 
$\ddt: \C Y \longrightarrow \C Y \otimes \C Y$. 
Since $\C Y$ is the set of generators of the algebras $\Hp$ and $\He$, 
and $\ddt(\|)=\| \otimes \|$, we can extend $\ddt$ to two maps 
$\ddp: \Hp \longrightarrow \Hp \otimes \Ha$ and 
$\dde: \He \longrightarrow \He \otimes \Ha$ defined as $\ddt$ 
on the generators (single trees), extended multiplicatively on tensor 
products, 
\begin{eqnarray*}
	\ddp (t_1 \cdots t_n) = \dde (t_1 \cdots t_n) 
	&:=& \ddt (t_1) \cdots \ddt (t_n),  
\end{eqnarray*}
and finally passed to the quotent $\Hat \longrightarrow \Ha$. 
Explicitly, $\ddp$ and $\dde$ can be recursively defined as 
\begin{eqnarray}
\label{ddpe} 
	\ddp (t \vee s) &=& \sum_{\Da t, \ddp s} 
	t_{(1)} \vee s_{(\gamma)} \otimes t_{(2)} \lpro s_{(\alpha)} , \\ 
	\dde (t \vee s) &=& \sum_{\Da t, \dde s} 
	t_{(1)} \vee s_{(e)} \otimes t_{(2)} \lpro s_{(\alpha)}, 
\end{eqnarray}
where we use the Sweedler notations 
\begin{eqnarray*} 
	\ddp s = \sum s_{(\gamma)} \otimes s_{(\alpha)}, &\quad& 
	\dde s = \sum s_{(e)} \otimes s_{(\alpha)}. 
\end{eqnarray*}
\end{ale}

\begin{ela}{\bf Lemma.} 
\label{lemmaddpe}
The maps $\ddp$ and $\dde$ are right $\Da$-coactions, i.e. they satisfy  
(\ref{D1-coassociative}), and they commute respectively with $\DPp$ and 
$\DPe$, i.e. they satisfy (\ref{D2-commute}). 
\end{ela}

\begin{proof}
Since the proof is exactly the same in the two cases, we do it 
explicitly only for $\dde$. 

The map $\dde$ is a right $\Da$-coaction, because we already proved 
that the identity 
\begin{eqnarray*}
	(\dde \otimes \Id) \dde &=& (\Id \otimes \Da) \dde 
\end{eqnarray*}
holds on single trees, and on a product $t_1\cdots t_n$, it follows from 
the fact that 
$$
	(\dde \otimes \Id) \dde(t_1\cdots t_n) = 
	[(\dde \otimes \Id) \dde(t_1)] \cdots [(\dde \otimes \Id) \dde(t_n)]
$$
and similarly 
$$
	(\Id \otimes \Da) \dde(t_1\cdots t_n) = 
	[(\Id \otimes \Da) \dde(t_1)] \cdots [(\Id \otimes \Da) \dde(t_n)]. 
$$

Let us prove that $\dde$ commutes with $\DPe$, i.e. that 
\begin{eqnarray*}
	(\DPe \otimes \Id) \dde &=& 
	m_{24}^3 (\dde \otimes \dde) \DPe , 
\end{eqnarray*}
where $m_{24}^3$ is the commutative multiplication in $\Ha$ with the 
notations of (\ref{semidirect Hopf}). 

On single trees, we prove it by induction. It is true for the root tree $\|$, 
so let us suppose that the equality holds for all trees up to order $n$, 
and let $t \vee s$ has order $n+1$. Then, on the left hand side we have 
\begin{eqnarray*}
	(\DPe \otimes \Id) \dde (t \vee s) 
	&=& \sum_{\Da t, \dde s} \DPe(t_{(1)} \vee s_{(e)}) 
	\otimes t_{(2)} \lpro s_{(\alpha)} \\ 
	&=& \sum_{\Da t, \dde s} \| \otimes t_{(1)} \vee s_{(e)} 
	\otimes t_{(2)} \lpro s_{(\alpha)} 
	+ \sum_{{\Da t, \dde s}\atop{\DPe s_{(e)}}} 
	t_{(1)} \vee {s_{(e1)}} 
	\otimes {s_{(e2)}} \otimes t_{(2)} \lpro s_{(\alpha)}, 
\end{eqnarray*}
while on the right hand side we have 
\begin{eqnarray*}
&& \hspace{-.5cm} m_{24}^3 (\dde \otimes \dde) \DPe (t \vee s) 
	\ = \ m_{24}^3 (\dde \otimes \dde) \left[ \| \otimes t\vee s 
	+ \sum_{\DPe s} t \vee s_{(1)} \otimes s_{(2)} \right] \\ 
&& \hspace{1cm} = \  m_{24}^3 \left( \sum_{\Da t, \dde s} \| \otimes \| 
	\otimes t_{(1)} \vee s_{(1)} \otimes t_{(2)} \lpro s_{(2)}  
	+ \sum_{\DPe s} \dde(t \vee s_{(1)}) \otimes \dde( s_{(2)}) \right) \\ 
&& \hspace{1cm} = \sum_{\Da t, \dde s} \| \otimes t_{(1)} \vee s_{(e)} 
	\otimes t_{(2)} \lpro s_{(\alpha)} 
	+ \sum_{{\DPe s,\Da t}\atop{\dde s_{(1)},\dde s_{(2)}}}
	t_{(1)} \vee {s_{(1e)}} \otimes s_{(2e)}
	\otimes t_{(2)} \lpro s_{(1\alpha)} \lpro s_{(2\alpha)} .  
\end{eqnarray*}
Then the two sides coincide, because for the tree $s$ we know that 
\begin{eqnarray*}
	\sum_{\dde s, \DPe s_{(e)}} 
	s_{(e1)} \otimes s_{(e2)} \otimes s_{(\alpha)} &=& 
	\sum_{\DPe s,\dde s_{(1)},\dde s_{(2)}} 
	s_{(1e)} \otimes s_{(2e)} \otimes s_{(1\alpha)} \lpro s_{(2\alpha)}. 
\end{eqnarray*}

Finally, we prove that the equality hold on a tensor product $t s \in \He$. 
On one side we have 
\begin{eqnarray*}
	(\DPe \otimes \Id) \dde (t s) 
	&=& (\DPe \otimes \Id) \left[ 
	\sum_{\dde t, \dde s} t_{(e)} s_{(e)}) 
	\otimes t_{(\alpha)} \lpro s_{(\alpha)} \right] \\ 
	&=& \sum_{{\dde t, \dde s}\atop{\DPe t_{(e)}, \DPe s_{(e)}}} 
	t_{(e1)} s_{(e1)} \otimes t_{(e2)} s_{(e2)} 
	\otimes t_{(\alpha)} \lpro s_{(\alpha)} \\ 
	&=& \left[ (\DPe \otimes \Id) \dde (t) \right] \ 
	\left[ (\DPe \otimes \Id) \dde (s) \right] .  
\end{eqnarray*} 
On the other side we have 
\begin{eqnarray*}
&& \hspace{-1.5cm} m_{24}^3 (\dde \otimes \dde) \DPe (t s) 
	\ =\ m_{24}^3 (\dde \otimes \dde) \left[ 
	\sum_{\DPe t,\DPe s} t_{(1)} s_{(1)}) \otimes t_{(2)} s_{(2)} \right]\\
&& \ = \ m_{24}^3 \left[ \sum_{{\DPe t, \DPe s}\atop
	{\dde t_{(1)}, \dde s_{(1)}, \dde t_{(2)}, \dde s_{(2)}}} 
	t_{(1e)} s_{(1e)}) \otimes t_{(1\alpha)} \lpro s_{(1\alpha)} 
	\otimes t_{(2e)} s_{(2e)} \otimes 
	t_{(2\alpha)} \lpro s_{(2\alpha)} \right] \\ 
&& \ = \ \sum_{{\DPe t, \DPe s}\atop
	{\dde t_{(1)}, \dde s_{(1)}, \dde t_{(2)}, \dde s_{(2)}}} 
	t_{(1e)} s_{(1e)} \otimes t_{(2e)} s_{(2e)} \otimes 
	t_{(1\alpha)} \lpro s_{(1\alpha)} 
	\lpro t_{(2\alpha)} \lpro s_{(2\alpha)}, 
\end{eqnarray*} 
which is equal to 
\begin{eqnarray*}
	&& \left[ m_{24}^3 (\dde \otimes \dde) \DPe (t) \right] \ 
	\left[ m_{24}^3 (\dde \otimes \dde) \DPe (s) \right]   
\end{eqnarray*} 
because $\lpro$ is commutative in $\Ha$. 
Then the two sides coincide by inductive hypothesis. 
\end{proof}

\begin{ale}{\bf The QED Hopf algebra.}
\label{HopfQED}
By Molnar's result of \cite{Molnar}, the smash coproduct 
$\Hq :=\Ha\ltimes\He$, as defined in (\ref{semidirect Hopf}), 
is then a graded connected Hopf algebra, which is neither commutative 
nor cocommutative.   
The grading is given by the sum of the orders of all the trees 
appearing in a monomial. 

The coproduct $\Dq: \Hq \longrightarrow \Hq \otimes \Hq$ is 
explicitly given by  
\begin{eqnarray*}
\label{Dq}
	\Dq(t \otimes s_1 \ldots s_n) &:=& 
	\Da(t) \ [(\dde \otimes \Id) \DPe (s_1 \ldots s_n)] .   
\end{eqnarray*}
\end{ale}

\begin{ale}{\bf The electron renormalization coaction.}
\label{coactionDe}
As in (\ref{semidirect Hopf}), we can then define a coaction 
of $\Hq$ on $\He$, as the map $\De: \He \longrightarrow \He \otimes \Hq$ 
given by 
\begin{eqnarray*}
\label{De}
	\De(s_1 \ldots s_n) &:=& 
	(\dde \otimes \Id) \DPe (s_1 \ldots s_n). 
\end{eqnarray*}
For instance, 
\begin{eqnarray*}
\De \| &=& \| \otimes \| \otimes \| \\ 
\De \Y &=& \Y \otimes \| \otimes \|+ \|\otimes \| \otimes \Y \\
\De \deuxun &=& 
    \deuxun \otimes \| \otimes \|
     + \Y \otimes \Y \otimes \| 
     + \| \otimes \| \otimes \deuxun \\
\De \deuxdeux &=& 
    \deuxdeux \otimes \| \otimes \|
     + \Y \otimes \| \otimes \Y 
     + \| \otimes \| \otimes \deuxdeux \\
\De\troisun &=& 
    \troisun \otimes \| \otimes \|
     + 2 \deuxun \otimes \Y \otimes \|
     + \Y \otimes \deuxun \otimes \| 
     + \| \otimes \| \otimes \troisun,\\
\De\troisdeux &=& 
    \troisdeux \otimes \| \otimes \|
     + \Y \otimes \deuxdeux \otimes \|
     + \| \otimes \| \otimes \troisdeux\\
\De\troistrois &=& 
    \troistrois \otimes \| \otimes \|
     + \deuxdeux \otimes \Y \otimes \|
     + \deuxun \otimes \| \otimes \Y
     + \Y \otimes \Y \otimes\Y 
     + \| \otimes \| \otimes \troistrois\\
\De\troisquatre &=& 
    \troisquatre \otimes \| \otimes \|
     + \deuxdeux \otimes \Y \otimes \|
     + \Y \otimes \| \otimes \deuxun 
     + \| \otimes \| \otimes \troisquatre,\\
\De\troiscinq &=& 
    \troiscinq \otimes \| \otimes \|
     + \deuxdeux \otimes \| \otimes \Y
     + \Y \otimes \| \otimes \deuxdeux
     + \| \otimes \| \otimes\troiscinq.
\end{eqnarray*}
\end{ale} 

\begin{ela}{\bf Lemma.} 
The coaction $\De$ of $\Hq$ on $\He$ can be defined recursively as 
\begin{eqnarray*}
	\De \| &=& \| \otimes \| \otimes \|, \\ 
	\De (t \vee s) &=& \| \otimes \| \otimes t \vee s 
	+ \sum_{\Da(t),\De(s)} t_{(1)} \vee s_{(1)} \otimes 
	t_{(2)} \lpro s_{(2)} \otimes s_{(3)} , \\ 
	\De (s_1 \ldots s_n) &=& \De (s_1) \cdots \De (s_n), 
\end{eqnarray*}
where we adopt a Sweedler's notation 
$\De s = \sum s_{(1)} \otimes s_{(2)} \otimes s_{(3)}$.  
\end{ela}
 
\begin{proof}
Since the maps $\DPe$ and $\dde$ are algebra morphisms, we only need to 
show it on the generators. We show it by induction on the order of trees. 
It is true for the root tree $\|$. Suppose that the equality holds for 
all trees up to order $n$, and let $t \vee s$ have order $n+1$. 
Then, in particular for $s$, we know that 
\begin{eqnarray*}
	\De s = \sum s_{(1)} \otimes s_{(2)} \otimes s_{(3)} &=&
	\sum_{\DPe(s), \dde(s_{(2)})} 
	s_{(1e)} \otimes s_{(1\alpha)} \otimes s_{(2)} . 
\end{eqnarray*}
So, applying the definition of $\De$ on $t \vee s$ and using the recursive 
definition (\ref{right pruning coproduct}) for $\DPe$ and 
(\ref{ddpe}) for $\dde$, we have 
\begin{eqnarray*}
	\De(t \vee s) &=& (\dde \otimes \Id) \DPe (t \vee s)  
	= \| \otimes \| \otimes t \vee s 
	+ \sum_{\DPe s} \dde(t \vee s_{(1)}) \otimes s_{(2)}\\ 
	&=& \| \otimes \| \otimes t \vee s + 
	\sum_{\DPe s,\Da t,\dde s_{(2)} } 
	t_{(1)} \vee s_{(1e)} \otimes t_{(2)} \lpro s_{(1\alpha)} 
	\otimes s_{(2)} \\ 
	&=& \| \otimes \| \otimes t \vee s 
	+ \sum_{\Da(t),\De(s)} t_{(1)} \vee s_{(1)} \otimes 
	t_{(2)} \lpro s_{(2)} \otimes s_{(3)}. 
\end{eqnarray*}
\end{proof}

\begin{ale}{\bf The photon renormalization coaction.} 
Exactly as in (\ref{HopfQED}), the semi-direct coproduct $\Ha \ltimes \Hp$ is 
a graded connected Hopf algebra, with twisted coproduct 
\begin{eqnarray*} 
	&& t \otimes s_1 \ldots s_n \mapsto 
	\Da(t) \ [(\dde \otimes \Id) \DPp (s_1 \ldots s_n)] , 
\end{eqnarray*}
which coacts on $\Hp$ from the right, with coaction given by the restriction 
of the coproduct to the subspace $\Hp$, as in (\ref{semidirect Hopf}) and 
(\ref{coactionDe}). 

However it is not the semidirect coproduct $\Ha \ltimes \Hp$ which describes 
the renormalization of the photon propagators. 
As we sketched in (\ref{renormalization group}), the photon renormalization 
Hopf algebra is the charge algebra $\Ha$, and the coaction is a semidirect 
coproduct coaction induced by a 1-cocycle, as in (\ref{semidirect Hopf}). 

Let $\sigma : \Hp \longrightarrow \Ha$ be the algebra morphism defined by 
\begin{eqnarray*}
	\sigma (t_1 \ldots t_n) &:=& t_1 \lpro\cdots\lpro t_n. 
\end{eqnarray*}
Then define $\Dp: \Hp \longrightarrow \Hp \otimes \Ha$ as the map 
\begin{eqnarray*}
	\Dp &:=& m_{23}^3 (\ddp \otimes \sigma) \DPp. 
\end{eqnarray*}
Since $\sigma$ is an algebra morphism, $\Dp$ is also an algebra morphism. 
\end{ale}

\begin{ela}{\bf Lemma.} 
The map $\Dp$ is a right coaction of $\Ha$ on $\Hp$, i.e. it is 
coassociative with respect to $\Da$. 
\end{ela}

\begin{proof} 
By lemma (\ref{sigma}), it is sufficient to show that the map $\sigma$ 
interwines $\Da$ and $\Dp$, i.e. for any $t_1 \ldots t_n \in \Hp$ we have 
\begin{eqnarray*}
	\Da \sigma (t_1 \ldots t_n) &=& 
	(\sigma \otimes \Id) \Dp(t_1 \ldots t_n). 
\end{eqnarray*}
If $n>1$, the result follows from the fact that all the maps are algebra 
morphisms.  So we only need to check it on a single tree $t$, for which 
$\sigma(t)=t$. 
We prove it by induction on the order of the trees. The equality holds 
for $\|$, suppose that it holds for a tree $t$, i.e. that 
\begin{eqnarray*}
	\sum_{\Da t} t_{(1)} \otimes t_{(2)} &=& 
	\sum_{\DPp t, \ddp t_{(1)}} \sigma(t_{(1 \gamma)})  
	\otimes t_{(2)} \lpro t_{(1\alpha)} 
	= \sum_{\DPp t, \ddp t_{(1)}} t_{(1 \gamma)} 
	\otimes t_{(2)} \lpro t_{(1\alpha)} . 
\end{eqnarray*}
Then for a tree $t \vee s$ with larger order we have: 
\begin{eqnarray*}
	(\sigma \otimes \Id) \Dp(t \vee s) &=& 
	(\sigma \otimes \Id) m_{23}^3 (\ddp \otimes \Id) \DPp (t \vee s) \\ 
	&=& (\sigma \otimes \Id) m_{23}^3 \left[ \ddp(t\vee s) \otimes \| 
	+ \sum_{\DPp t} \ddp(t_{(1)}) \otimes t_{(2)} \vee s \right] \\ 
	&\stackrel{(1)}{=}& \ddp(t\vee s) + \sum_{\DPp t, \ddp t_{(1)}} 
	\sigma(t_{(1\gamma)}) \otimes t_{(1\alpha)} \lpro (t_{(2)} \vee s) \\ 
	&=& \ddp(t\vee s) + \sum_{\Da t} 
	t_{(1)} \otimes t_{(2)} \lpro V(s) = \Da (t \vee s) , 
\end{eqnarray*}
where the equality $(1)$ holds because $\ddp$ applied to a single tree 
produces only single tree components on the left hand side, hence 
$(\sigma \otimes \Id) \ddp(t\vee s) = \ddp(t\vee s)$. 
\end{proof} 

Remark that the coaction $\Dp$ applied to a single tree also 
produces only single tree components on the left hand side, hence 
$(\sigma \otimes \Id) \Dp(t) = \Dp(t)$ for any $t \in \C Y$. 
In conclusion, we obtain a very simple expression for the 
photon renormalization coaction on single trees. 

\begin{ela}{\bf Corollary.}
The photon renormalization coaction $\Dp$ restricted to the 
subspace of single trees coincides with the non-commutative 
charge coproduct $\Dat$, 
\begin{eqnarray*}
	\Dp(t) &=& \Dat(t), \qquad\mbox{for any $t \in Y$}. 
\end{eqnarray*}
\end{ela}

Exemples of $\Dp(t)$ for small order trees can then be constructed 
directly from (\ref{HopfHa}): 
\begin{eqnarray*}
	\Dp \| &=& \| \otimes \| \\ 
	\Dp \Y &=& \Y\otimes \| + \| \otimes \Y , \\ 
	\Dp \deuxun &=& \deuxun\otimes \| + 2\Y \otimes \Y 
	+ \| \otimes \deuxun , \\ 
	\Dp \deuxdeux &=& \deuxdeux\otimes \| + \| \otimes \deuxdeux , \\ 
	\Dp\troisun &=& \troisun \otimes \| + 3 \deuxun \otimes \Y 
	+ 3 \Y\otimes \deuxun + \| \otimes \troisun ,\\
	\Dp\troisdeux &=& \troisdeux \otimes \| + \deuxdeux\otimes \Y 
	+ \Y\otimes \deuxdeux + \| \otimes \troisdeux ,\\
	\Dp\troistrois &=& \troistrois\otimes \| + \deuxdeux\otimes \Y 
	+ \Y \otimes \deuxdeux + \| \otimes \troistrois ,\\
	\Dp\troisquatre &=& \troisquatre \otimes \| + \deuxdeux \otimes \Y 
	+ \| \otimes \troisquatre ,\\
	\Dp\troiscinq &=& \troiscinq \otimes \| + \| \otimes \troiscinq .
\end{eqnarray*}


\section{Renormalization of tree-expanded QED propagators} 
\label{ren-tree}

Let $\alpha_0$ be the bare fine structure constant (before renormalization).  
For each momentum vector $q \in C^4$, let $D(\alpha_0;q)$ and $S(\alpha_0;q)$ 
denote the bare Feynman propagators for the photon and electron fields, 
as consi\-dered in \cite{BrouderFrabetti2}\footnote{
Remark that $D(\alpha_0;q)$ and $S(\alpha_0;q)$ differ from the usual QED 
propagators respectively by a factor $D_0(q)^{-1}$ and $S_0(q)^{-1}$, 
which are the inverse of the free propagators.}. 
Following \cite{Brouder}, consider the {\em tree-expansions}  
\begin{eqnarray}
	D(\alpha_0;q) &=& \sum_{t \in Y} \Up_q(t) \alpha_0^{|t|} , 
	\label{D} \\ 
	S(\alpha_0;q) &=& \sum_{t \in Y} \Ue_q(t) \alpha_0^{|t|} , 
	\label{S}  
\end{eqnarray} 
that is, the expansions of these propagators as power series of $\alpha_0$ 
with coefficients labeled by planar binary trees. 
For single particles, the coefficients $\Up_q(t)$ and $\Ue_q(t)$ are 
$4 \times 4$ complex matrices. The coefficents $\Up_q(\|)$ and $\Ue(\|)$ 
of the root tree represent the free propagators, which, by assumption, 
are the identity $4 \times 4$ matrix $I$\footnote{
Because $\Up_q(\|) = D_0(q) D_0(q)^{-1} = I$ 
and similarly $\Ue_q(\|) = S_0(q) S_0(q)^{-1} = I$.}. 
For higer order trees, the coefficients $\Up_q(t)$ and $\Ue_q(t)$ 
can be explicitly determined as Feynman amplitudes, since each tree 
is a finite sum of appropriate Feynman diagrams, cf. \cite{BrouderFrabetti2}. 
In alternative, they can be determined recursively, as showed in 
\cite{Brouder}, starting from the coefficients of the smaller trees 
$t^l$ and $t^r$ such that $t^l \vee t^r=t$. 
In conclusion, the tree-expansions (\ref{D},\ref{S}) allow to consider 
trees as a basis for some ``polynomial functions'' on the set of propagators, 
and therefore to identify the QED propagators with two algebra morphisms 
$\Up_q$ and $\Ue_q$, on $\Hp$ and $\He$ respectively, such that 
\begin{eqnarray*}
	\< t,D(\alpha_0;q) \> &\equiv& \< \Up_q,t \> = \Up_q(t), \\ 
	\< t,S(\alpha_0;q) \> &\equiv& \< \Ue_q,t \> = \Ue_q(t) . 
\end{eqnarray*}
Moreover, in \cite{Brouder} it was shown that the product of propagators 
is dual to the pruning coproducts, that is 
\begin{eqnarray*}
	\< t,D(\alpha_0;q) D(\alpha_0;q) \> &=& 
	\< \Up_q \otimes \Up_{q},\DPp (t) \>, \\ 
	\< t,S(\alpha_0;q) S(\alpha_0;q) \> &=& 
	\< \Ue_q \otimes \Ue_{q}, \DPe(t) \> . 
\end{eqnarray*}

Furthermore, let $\alpha$ be the renormalized fine structure constant, and 
let $\bar D(\alpha;q)$ and $\bar S(\alpha;q)$ denote the 
massless renormalized propagators as in \cite{BrouderFrabetti2}. 
Again, consider the tree-expansions 
\begin{eqnarray}
	\bar D(\alpha;q) &=& \sum_{t \in Y} \Rp_q(t) \alpha^{|t|} 
	\label{barD} \\ 
	\bar S(\alpha;q) &=& \sum_{t \in Y} \Re_q(t) \alpha^{|t|} 
	\label{barS}  
\end{eqnarray} 
as power series on $\alpha$, also starting with the unperturbed 
coefficient given by $I$. As before, these expansions determine two 
algebra morphisms $\Rp_q$ and $\Re_q$ on $\Hp$ and $\He$ respectively. 
The aim of renormalization theory is to find their values $\Rp_q(t)$ 
and $\Re_q(t)$ on all the trees. 
In \cite{BrouderFrabetti1}, we gave some recursive solutions\footnote{
The recursive solutions given in \cite{BrouderFrabetti1} are valid 
in massive renormalization.} with respect to the order of the trees. 
Here we recall how the relationship between all these coefficients can be 
given in terms of the Hopf algebras and coactions on trees defined 
in the previous sections. 

Let $Z_3(\alpha)$ and $Z_2(\alpha)$ denote the renormalization 
factors for the photon and the electron propagators. They satisfy 
the Dyson formulas 
\begin{eqnarray*}
	\bar D(\alpha;q) Z_3(\alpha) &=& D(\alpha_0;q) \\ 
	\bar S(\alpha;q) Z_2(\alpha) &=& S(\alpha_0;q) 
\end{eqnarray*} 
and the charge renormalization formula proved by Ward 
\begin{eqnarray*}
	\alpha_0(\alpha) &=& \alpha Z_3(\alpha)^{-1}. 
\end{eqnarray*} 

As explained in \cite{BrouderFrabetti1}, \cite{BrouderFrabetti2}, 
trees represent sums of Feynman diagrams. It is well known that 
the renormalization factors are expanded only over 1PI Feynman graphs, 
and this property corresponds to the following expansions over trees: 
\begin{eqnarray*}
	Z_3(\alpha) &=& 1 - \sum_{t \in Y} \Cp(V(t)) \alpha^{|t|+1}, \\ 
	Z_2(\alpha) &=& 1 + \sum_{t \neq \|} \Ce(\SPe(t)) \alpha^{|t|},  
\end{eqnarray*} 
where $V(t)$ are the generators of the algebra $\Ha$ and 
$\SPe(t)$ are the elements of the algebra $\He$ transformed under 
the right pruning antipode defined in (\ref{pruning Hopf}). 
Once again, these expansions determine two algebra morphisms $\Cp$ 
and $\Ce$ on $\Hp$ and $\He$ respectively. In the present case they are 
both scalars, but the formalism based on planar binary trees allows 
to consider also non-scalar maps. 
Moreover, the Ward formula for the fine structure constant tells us 
that the map $\Cp$ is also an algebra morphism on $\Ha$. 

We finally state a results of \cite{BrouderFrabetti1}, which shows 
that the coproduct on trees previously defined encodes the  
relationship between the amplitudes, before and after the renormalization. 

\begin{ela}{\bf Theorem.}
\label{theoremren}
The relation between the coefficients of the expansions 
(\ref{D}) and (\ref{barD}) for the bare and the renormalized 
photon propagators is
\begin{eqnarray*}
	\Rp_q(t) &=& \< \Up \otimes \Cp, \Dp(t) \> 
	= \sum_{\Dp(t)} \Up_q(t_{(1)}) \Cp(t_{(2)}) .  
\end{eqnarray*}
The relation between the coefficients of the expansions 
(\ref{S}) and (\ref{barS}) for the bare and the renormalized 
electron propagators is 
\begin{eqnarray*}
	\Re_q(t) &=& \< \De \otimes \Cp \otimes \Ce, \De(t) \>
	= \sum_{\De(t)} \Ue(t_{(1)};q) \Cp(t_{(2)}) \Ce(\SPe t_{(3)}) .
\end{eqnarray*}
\end{ela}

It remains to show that the Hopf algebra $\Ha$ describes the renormalization 
of the fine structure constant $\alpha$, i.e. that 
\begin{eqnarray*}
	\< (\alpha_1 \circ \alpha_2)(\alpha), t \> 
	&=& \< \alpha_1 \otimes \alpha_2 , \Da t \>,  
\end{eqnarray*}
if $\alpha_2(\alpha)$ and $\alpha_1(\alpha_2)$ are two successive 
renormalizations. This is the topic of the paper \cite{BrouderFrabetti4} 
in preparation, where we define the natural group of composition of 
series expanded over trees. 


\section{Conclusions} 
\label{conclusions}

In \cite{ConnesKreimerII}, Connes and Kreimer relate the renormalization 
group of the $\Phi^3$ theory, based on Feynman graphs, to the group of 
formal diffeomorphisms on the complex line, and to the Birkhoff 
decomposition of holomorphic line bundles on the circle. 
In their case, since the renormaliation Hopf algebra is commutative, 
the Milnor-Moore theorem allows to study the dual renormalization group 
even without knowing it explicitly. In our case, since the renormalization 
Hopf algebra of QED propagators is neither commutative nor cocommutative, 
we need first to introduce its natural dual group, which {\em a priori} 
does not necessarily exist. This is the topic of the paper 
\cite{BrouderFrabetti4} in preparation. 

However, the QED Hopf algebra on trees can be directly related to the Hopf 
algebra dual to the group of formal diffeomorphisms. This comparison 
can be done simply by summing up all the trees at a given order $n$, 
which corresponds to the order $n$ of interaction for the particle Green 
functions. The result is a non-commutative version of the Hopf algebra 
of formal diffeomorphisms, described in \cite{BrouderFrabetti3}. 

Finally, all these results suggest a natural question which has no answer 
yet: can the renormalization group of perturbative quantum field theories 
be realized as an ``automorphism group'' on some space?  



\begin{thebibliography}{10}
\addcontentsline{toc}{section}{\bf References}

\bibitem{Brouder}
Ch. Brouder. 
\newblock {\em On the trees of quantum fields}, 
\newblock Eur. Phy. J. C, {\bf 12} (2000), 535-549.

\bibitem{BrouderFrabetti1}
Ch. Brouder and A.~Frabetti.
\newblock {\em Renormalization of {QED} with planar binary trees}, 
\newblock Eur. Phys. J. C {\bf 19} (2001), 715-741.

\bibitem{BrouderFrabetti2}
Ch. Brouder and A.~Frabetti.
\newblock {\em Noncommutative renormalization of massless {QED}}, 
\newblock hep-th/0011161. 

\bibitem{BrouderFrabetti3}
Ch. Brouder and A.~Frabetti.
\newblock {\em Noncommutative Hopf algbera of formal diffeomorphisms}, 
\newblock in preparation. 

\bibitem{BrouderFrabetti4}
Ch. Brouder and A.~Frabetti.
\newblock {\em Groups of tree-expanded series}, 
\newblock in preparation. 

\bibitem{Butcher}
J.C.~Butcher.  
\newblock {\em An algebraic theory of integration methods}, 
\newblock Math. Comput., {\bf 26} (1972), 79-106. 

\bibitem{ConnesKreimerIHES}
A.~Connes and D.~Kreimer. 
\newblock {\em Hopf algebras, renormalization and noncommutative geometry}, 
\newblock Comm. Math. Phys. {\bf 199} n.1 (1998), 203-242.

\bibitem{ConnesKreimerI}
A.~Connes and D.~Kreimer. 
\newblock {\em Renormalization in quantum field theory and the 
{R}iemann-{H}ilbert problem {I}: the {H}opf algebra structure 
of graphs and the main theorem}, 
\newblock Comm. Math. Phys. {\bf 210} n.1 (2000), 249-273.

\bibitem{ConnesKreimerII}
A.~Connes and D.~Kreimer.
\newblock {\em Renormalization in quantum field theory and the 
{R}iemann-{H}ilbert problem {II}: the $\beta$ function, diffeomorphisms 
and the renormalization group}, 
\newblock Comm. Math. Phys. {\bf 216} n.1 (2001), 215-241. 

\bibitem{Foissy}
L.~Foissy. 
\newblock {\em Les alg\`ebres de Hopf des arbres enracin\'es
decor\'es}, 
\newblock Th\`ese Univ. Reims (2001), math.QA/0105212. 

\bibitem{Frabetti}
A. Frabetti. 
\newblock {\em Simplicial properties of the set of planar binary trees},  
\newblock J. Alg. Comb. {\bf 13} (2001), 41-65. 

\bibitem{GrossmanLarson}
R.~Grossman and R.G.~Larson. 
\newblock {\em Hopf algebraic structure of families of trees}, 
\newblock J. Algebra {\bf 126} n.1 (1989), 184-210. 

\bibitem{Holtkamp}
R.~Holtkamp.  
\newblock {\em Comparison of Hopf algebras on trees},  
\newblock Preprint (2001). 

\bibitem{Kreimer98}
D.~Kreimer.
\newblock {\em On the {H}opf algebra structure of perturbative quantum 
field theory}, 
\newblock Adv. Th. Math. Phys. {\bf 2} (1998), 303-334.

\bibitem{ItzyksonZuber}
C. Itzykson and J.-B. Zuber.
\newblock {\em Quantum Field Theory}, 
\newblock McGraw-Hill, New York, 1980. 

\bibitem{Lin}
B.~Lin. 
\newblock {\em Crossed coproducts of Hopf algebras}, 
\newblock Comm. Alg. {\bf 10} n.1 (1982), 1-17. 

\bibitem{LodayATree}
J.-L. Loday.
\newblock {\em Arithmetree}, 
\newblock Preprint math.CO/0112034

\bibitem{LodayDias}
J.-L. Loday. 
\newblock {\em Alg\`ebres ayant deux operations associatives (dig\`ebres)},  
\newblock C. R. Acad. Sci. Paris {\bf 321} (1995), 141-146. 

\bibitem{LodayRonco1}
J.-L. Loday and M.O. Ronco.
\newblock {\em Hopf algebra of the planar binary trees}, 
\newblock Adv. Math. {\bf 139} (1998), 293-309.

\bibitem{LodayRonco2}
J.-L. Loday and M.O. Ronco.
\newblock {\em Order structure and the algebra of permutations and 
planar binary trees}, 
\newblock J. Alg. Comb., to appear. 

\bibitem{Majid}
S.~Majid. 
\newblock {\em Physiscs for algebraists: noncommutative and 
noncocommutative Hopf algebras by a bicrossedproduct}, 
\newblock J. Alg. {\bf 130} n.1 (1990), 17-64. 

\bibitem{Moerdijk}
I.~Moerdijk.  
\newblock {\em On the Connes-Kreimer construction of Hopf algebras}, 
\newblock Cont. Math. {\bf 271} (2001), 311-321. 

\bibitem{Molnar}
R.K.~Molnar. 
\newblock {\em Semi-direct products of Hopf algebras}, 
\newblock J. Alg. {\bf 47} (1977), 29-51. 

\bibitem{Panaite}
F.~Panaite. 
\newblock {Relating the Connes-Kreimer and the Grassman-Larson Hopf
algebras on rooted trees}, 
\newblock Lett. Math. Phys. {bf 51} n.3 (2000), 211-219. 

\bibitem{PeskinSchroeder}
M.E.~Peskin and D.V.~Schroeder. 
\newblock {\em An Introduction to Quantum Field Theory}, 
\newblock Perseus Books Pub. L.L.C., 1995. 

\bibitem{Radford}
D.~Radford. 
\newblock {\em The structure of Hopf algebras with a projection}, 
\newblock J. Alg. {\bf 92} n.2 (1985), 322-347. 

\bibitem{Ward}
J.C.~Ward. 
\newblock {\em An identity in quantum electrodynamics}, 
\newblock Phys. Rev. {\bf 78} (1950), 182. 

\end{thebibliography}
\end{document}